\magnification1200
\input amssym.def
\input amssym.tex 
\def\SetAuthorHead#1{}
\def\SetTitleHead#1{}
\def\NoindentAfter{\everypar={\setbox0=\lastbox\everypar={}}}
\def\H#1\par#2\par{{\baselineskip=15pt\parindent=0pt\parskip=0pt
 \leftskip= 0pt plus.2\hsize\rightskip=0pt plus.2\hsize
 \bf#1\unskip\break\vskip 4pt\rm#2\unskip\break\hrule
 \vskip40pt plus4pt minus4pt}\NoindentAfter}
\def\HH#1\par{{\bigbreak~\bigbreak\noindent\bf#1\medskip}\NoindentAfter}
\def\HHH#1\par{{\bigbreak\noindent\bf#1\unskip.\kern.4em}}
\def\th#1\par{\medbreak\noindent{\bf#1\unskip.\kern.4em}\it}

\def\pf#1\par{\medbreak\noindent{\it#1\unskip.\kern.4em}}
\def\df#1\par{\medbreak\noindent{\it#1\unskip.\kern.4em}}

\def\skipaline{\medskip}

\let\Roster\bgroup\def\endRoster{\egroup\par}
\def\\{}\def\text#1{\hbox{\rm #1}}

\def\MaxReferenceTag#1{}
\def\qedbox{\vrule width2mm height2mm\hglue1mm\relax}
\def\qed{\ifmmode\qedbox\else\hglue5mm\unskip\hfill\qedbox\medbreak\fi\rm}

\def\cite#1{{\bf[#1]}}

\def\Bib#1\par{\bigbreak\bgroup\centerline{#1}\medbreak\parindent30pt
 \parskip2pt\frenchspacing\par}
\def\endBib{\par\egroup}
\newdimen\Overhang
\def\rf#1{\par\noindent\hangafter1\hangindent=1 true in
     \setbox0=\hbox{[#1]}\Overhang\wd0\advance\Overhang.4em\relax
     \ifdim\Overhang>\hangindent\else\Overhang\hangindent\fi
     \hbox to \Overhang{\box0\hss}\ignorespaces}

\def\bbZ{{\Bbb Z}}
\newcount\notenumber\notenumber0
\def\note#1{\advance\notenumber1\footnote{$^{\the\notenumber)}$}{#1}}
\def\Coordinates{\bigbreak\bgroup\parindent=0pt\obeylines}
\def\endCoordinates{\egroup}
\overfullrule=0pt

\SetAuthorHead{Walter D. Neumann and Michael Shapiro}
\SetTitleHead{Regular geodesic normal forms}
\H
Regular geodesic normal forms in virtually abelian groups
\footnote{}{\rm AMS Classification 20F32}

Walter D. Neumann and Michael Shapiro

\footnote{}{Both authors acknowledge
support from the ARC for this research.}

\HHH Abstract

We describe a virtually abelian group $G$ generated by a finite set
$X$ such that there is no regular language of geodesic $X$-words that
surjects to $G$ by evaluation.

\HH 1.~~Introduction

Cannon gave an example in \cite{ECHLPT} (Example 4.4.2, see also
\cite{NS}) of a virtually abelian group with finite generating set $X$
such that the language of geodesic $X$-words for $G$ is not a regular
language.  Since his example admits a geodesic automatic structure in
the generators $X$, it left open the possibility that this might
always be so: a virtually abelian group might admit a geodesic
automatic structure for any generating set, or at least a geodesic
regular unique normal form (which is weaker).  This is suggested as a
question in \cite{ECHLPT}.  A reason to hope it might be true was that
it would give a very satisfactory proof of Benson's theorem \cite{B}
that the growth function of a virtually abelian group is rational with
respect to any generating set.

In this note we exhibit a virtually abelian group $G$ with finite
generating set $X$ such that there is {\it no} regular language of
geodesic $X$-words for $G$ that surjects to $G$.

In this context it is worth
recalling a result of \cite{NS} that any finite generating set of an
abelian group has the property that the full language of geodesics is
regular and contains an automatic structure for $G$.  Moreover, for
virtually abelian groups every finite generating set can be enlarged
to one with this property.
\skipaline
{\it Convention}. Our generating sets $X$ will always be symmetric.
That is, when we list generators we mean that both those generators
and their inverses are to be in our generating set.

\HH 2.~~A preliminary example

We start with an example that uses a weighted generating set.  That
is, we allow the generators to have lengths other than $1$.  Our
example is given by the following presentation:
$$\eqalign{H = \langle x,y,t,\tau \mid ~ & t^{2}=\tau^{2}=[\tau,t]=1,
\cr
& x,x^{t},x^{\tau},x^{t\tau},y,y^{t},y^{\tau},y^{t\tau} \text{
commute},\cr
&x+x^{\tau} = y + y^{t} \rangle. }$$
The last relation gives
$$\eqalign{ x+x^{\tau} &= y + y^{t}  \cr
            x+x^{\tau} &= y^{\tau} + y^{t\tau}\cr
            x^{t}+x^{t\tau} &= y + y^{t} \cr
            x^{t}+x^{t\tau} &= y^{\tau} + y^{t\tau}.}$$
These are not independent.  It is not hard to see that
$$1 \to N \to H \to V \to 1$$
where $V$ is the Klein 4-group and $N$ is free abelian of rank 5.
We will give the generators $x^{\pm1},y^{\pm1},t^{\pm1},\tau^{\pm1}$ 
lengths $1,1,1,2$
respectively.  We will suppose that there is a regular language $L$
consisting of geodesics and surjecting to $H$ and we will derive a
contradiction.

There is a homomorphism $\epsilon\colon H\to \bbZ$ that takes $t$ and
$\tau$ to the trivial element and
$x,x^t,x^\tau,x^{t\tau},y,y^t,y^\tau,y^{t\tau}$ all to $1$.  Any
element $h\in H$ thus has length at least $|\epsilon(h)|$.  Moreover,
if it is in $N$ then it involves an even number of instances of $t$
and of $\tau$, so its length is congruent to $\epsilon(h)$ modulo $2$.
Moreover, if it is in $N$ but not in the subgroup generated by $x$ and
$y$ then it has length at least $|\epsilon(h)|+2$, or even
$|\epsilon(h)|+4$ if it cannot be written in terms of $x$, $y$, and $t$.
 
We consider elements of $N<H$ in the positive span of $x$,
$x^{\tau}$, $y$ and $y^{t}$.  First consider an element of the form
$x^{a}y^{b}(y^{t})^{c}$ with $c>0$.  This can be written geodesically as
$x^{a}y^{b}ty^{c}t$ and thus $\ell(x^{a}y^{b}(y^{t})^{c})=a+b+c+2$.

We now consider an element of the form $x^{d}(x^{\tau})^{e}y^{f}$ with
$e>0$.  
Notice that
$$x^{d}(x^{\tau})^{e}y^{f} =
x^{d-\delta}(x^{\tau})^{e-\delta}y^{f+\delta}(y^{t})^{\delta}.$$ Thus,
if $e \le d$, this lies in the positive span of $x$, $y$, and $y^{t}$
and thus $\ell(x^{d}(x^{\tau})^{e}y^{f})=d+e+f+2$.  On the other hand,
if $e>d$, then it is not hard to see that $x^{d}(x^{\tau})^{e}y^{f}$
has length $d+e+f+4$.  And moreover, if $e>d+1$ then every geodesic for
this element has the form
$$w= w_{1} \tau x^{e}\tau w_{2}\eqno{(*)}$$ where $w_{1}$ and $w_{2}$
are positive words in $x$ and $y$ which together contain $d$ $x$'s
and $f$ $y$'s. (If $e=d+1$ we have additional geodesics such as
$x^{-1}y^{f+e}ty^et$ of length $2e+f+3=d+e+f+4$.)  Notice that the
word $(*)$ is geodesic if and only if $e>d$.  Now $L$ must contain
such words with $d$ arbitrarily large.  We fix a finite state automaton 
$M$ for $L$ and
choose a word $w$ so that $d$ is at least the number of states of $M$.
Then $e$ is greater than the number of states of $M$ so the portion
$x^e$ of $w$ traverses a loop of $M$.  By eliminating such loops we find
a word $w'=w_1\tau x^{e'}\tau w_2$ which is also in $L$, but has $e'$
no larger than the number of states of $M$ and hence no larger than $d$. It
is thus not geodesic.  Thus there can be no regular language of
geodesics surjecting to $H$.

\HH 3.~~The main example

We now embed the above preliminary example in an example where all
generators have length $1$. We take
$$\eqalign{G = \langle x,y,t,s \mid ~ & t^{2}=s^{4}=[s,t]=1,
\cr
& x^{\alpha},x^{\beta}, y^{\alpha},y^{\beta} \text{
commute, for any }\cr
&\alpha,\beta \in\{1,s,s^{2},s^{3},t,ts,ts^{2},ts^{3}\},\cr
&x+x^{s^{2}} = y + y^{t} \rangle . }$$
$G$ has a finite index normal subgroup which is free abelian of rank 10.
Taking $\tau=s^2$ embeds $H$ in to $G$.  It is easy to see that this
embedding is totally geodesic in the sense that a geodesic word in
$x,y,s,t$ that evaluates into the subgroup is the result of
substituting $\tau=s^2$ in a geodesic word of the subgroup.  We leave
this to the reader.  It then follows that this example inherits the
property that the language of geodesics has no regular sublanguage
which surjects to the group (in fact, no regular sublanguage can
evaluate to a subset containing the subgroup).

\Bib        Bibliography

\rf{B} M. Benson, Growth series of finite extensions of $\bbZ^n$ are
rational, Invent.\ Math.\ 73 (1983) 251--269.

\rf {ECHLPT} D.~B.~A. Epstein, J.~W. Cannon, D.~F. Holt, S.~V.~F. Levy,
M.~S. Paterson, and  W.~P. Thurston, Word processing in groups, Jones and
Bartlett, 1992.

\rf {NS} W.~D.~Neumann and M.~Shapiro, Automatic structures, rational growth,
and geometrically finite groups.  Invent. Math. {\bf120} (1995), 259--287.

\endBib
\Coordinates
Department of Mathematics\\
The University of Melbourne\\
Parkville, Vic 3052\\
Australia\\

Email: neumann@maths.mu.oz.au
\endCoordinates
\Coordinates
Department of Mathematics\\
The University of Melbourne\\
Parkville, Vic 3052\\
Australia\\

Email: shapiro@maths.mu.oz.au
\endCoordinates

\bye